\documentclass[leqno,12pt]{article}

\usepackage{latexsym}
 \usepackage{graphicx}
\usepackage{amsmath}
\usepackage{amssymb}
\usepackage{algorithm}
\usepackage{algorithmic}
\usepackage{cite}
\usepackage{color}

\usepackage[margin=1.3in, top=1.3in, bottom=1.3in]{geometry}

\newcommand{\nc}{\newcommand}
\nc{\nt}{\newtheorem}
\nt{thm}{Theorem}[section]
\nt{cor}[thm]{Corollary}
\nt{prop}[thm]{Proposition}
\nt{lem}[thm]{Lemma}
\nt{defn}[thm]{Definition}
\nt{rem}[thm]{Remark}
\nt{exa}[thm]{Example}
\nt{ass}[thm]{Assumption}
\nt{alg}[thm]{Algorithm}
\nt{que}[thm]{Question}
\nt{con}[thm]{Conjecture}
\nc{\ip}[2]{\mbox{$\langle #1,#2 \rangle$}}
\nc{\pf}{\noindent{\bf Proof\ \ }}
\nc{\finpf}{\hfill{$\Box$}\linespace}
\nc{\linespace}{\vspace{\baselineskip} \noindent}
\nc{\R}{{\bf R}}
\nc{\cM}{{\mathcal M}}
\nc{\X}{{\mathbf X}}
\nc{\Y}{{\mathbf Y}}
\nc{\Rn}{{\bf R}^n}
\nc{\Sn}{{\bf S}^n}
\nc{\lm}{\lambda_{\mbox{\rm\scriptsize min}}}
\nc{\bx}{\bar{x}}
\nc{\e}{\epsilon}
\nc{\cl}{\mbox{\rm cl}\,}
\nc{\conv}{\mbox{\rm conv}\,}

\def\tto{\;{\lower 1pt \hbox{$\rightarrow$}}\kern -12pt
           \hbox{\raise 2.8pt \hbox{$\rightarrow$}}\;}
\newenvironment{myequation}{\setcounter{equation}{\value{thm}}
   \begin{equation}}{\addtocounter{thm}{1}\end{equation}}

\nc{\bmye}{\begin{myequation}}
\nc{\emye}{\end{myequation}}

\begin{document}
\title{
Inexact alternating projections on nonconvex sets
}
\author{
 D. Drusvyatskiy
\thanks{Department of Mathematics, University of Washington, Seattle WA 98195, USA.
\texttt{ddrusv@uw.edu.}
Research supported in part by
AFOSR YIP award FA9550-15-1-0237 and by National Science Foundation Grants DMS 1651851 and CCF 1740551.
}
\and
A.S. Lewis
\thanks{ORIE, Cornell University, Ithaca, NY 14853, USA.
\texttt{people.orie.cornell.edu/aslewis}.
Research supported in part by National Science Foundation Grant DMS-1613996.}
}
\maketitle

\vspace{-\baselineskip}
\begin{center}
{\em
Dedicated to our friend, colleague, and inspiration, Alex Ioffe, \\
on the occasion of his 80th birthday.
}
\end{center}
\bigskip

\begin{abstract}
Given two arbitrary closed sets in Euclidean space, a simple transversality condition guarantees that the method of alternating projections converges locally, at linear rate, to a point in the intersection.  Exact projection onto nonconvex sets is typically intractable, but we show that computationally-cheap inexact projections may suffice instead.  In particular, if one set is defined by sufficiently regular smooth constraints, then projecting onto the approximation obtained by linearizing those constraints around the current iterate suffices.  On the other hand, if one set is a smooth manifold represented through local coordinates, then the approximate projection resulting from linearizing the coordinate system around the preceding iterate on the manifold also suffices.
\end{abstract}
\medskip

\noindent{\bf Key words:} Alternating projections, Linear convergence, Variational analysis, Transversality, Approximate projection
\medskip

\noindent{\bf AMS 2000 Subject Classification:} 
49M20, 65K10, 90C30
\medskip

\section{Nonconvex alternating projections}
The method of alternating projections for finding a point in the intersection of two closed convex sets enjoys an enduring popularity.  A suitable regularity or transversality condition ensures local linear convergence, even for nonconvex sets.  With a few exceptions (notably including manifolds of fixed-rank matrices), exact projection---finding exact nearest points---for nonconvex sets is practically unrealistic.  Nonetheless, the alternating projection philosophy remains appealing because we can often approximate the projection cheaply without losing linear convergence.

We present two examples, both feasibility problems involving a closed set $Q$, onto which we suppose we can project easily.  The first example concerns the constraint system
\[
\left\{
\begin{array}{rcl}
G(x)	& \le	& 0 \\
H(x)	& = 	& 0 \\
x		& \in	& Q,
\end{array}
\right.
\]
for $C^{2}$-smooth maps $G$ (into $\Rn$) and $H$.  For this problem, we consider the following algorithm, which involves just projections onto $Q$ and various polyhedra.

\begin{alg}
{\rm
\begin{algorithmic}
\STATE
\WHILE{$x \in Q$ not feasible}
\STATE  choose minimal norm $s$ satisfying\\
$\mbox{} \qquad G(x) + \nabla G(x)s \le 0$ and $H(x) + \nabla H(x)s = 0$;
\STATE  choose $x_+ \in Q$ minimizing $|x+s-x_+|$;
\STATE  $x=x_+$;
\ENDWHILE
\end{algorithmic}
}
\end{alg}

\noindent
The second example concerns the system
\[
F(x) \in Q,
\]
for a $C^{2}$-smooth map $F$.  Now we consider the following algorithm.

\begin{alg}
{\rm
\begin{algorithmic}
\STATE
\WHILE{$x$ not feasible}
\STATE	choose $y \in Q$ minimizing $|y-F(x)|$;
\STATE  choose $s$ minimizing $|F(x) + \nabla F(x)s - y|$;
\STATE	$x=x+s$;
\ENDWHILE
\end{algorithmic}
}
\end{alg}

\noindent
By interpreting each algorithm as an approximate alternating projection scheme, we derive reasonable conditions for local linear convergence to a feasible solution.

\section{The geometry of alternating projections} \label{alternating section}
We now describe our setting more formally.  Consider a Euclidean space ${\mathbf X}$ with norm $|\cdot|$.  
Given a nonempty closed set $Q \subset {\mathbf X}$, we define the {\em distance function} $d_Q \colon {\mathbf X} \to \R$ by
\[
d_Q(z) ~=~ \min \big\{ |x-z| : x \in Q \big\},
\]
for any point $z \in {\mathbf X}$.  The {\em projection} $P_Q(z)$ consists of the attaining (or {\em nearest}) points $x$:  it is a nonempty closed set, but may not be a singleton if $Q$ is nonconvex.  Given a second nonempty closed set $M \subset {\mathbf X}$, the {\em method of alternating projections} is the following iteration.

\begin{alg}[Alternating projections]
\label{alternating projections}
{\rm
\begin{algorithmic}
\STATE
\STATE  choose $z \in Q$; 
\WHILE{$z \not\in M$}
\STATE  choose $x \in P_M(z)$;
\STATE  choose $z' \in P_{Q}(x)$;
\STATE  $z=z'$;
\ENDWHILE
\end{algorithmic}
}
\end{alg}

\noindent
Consider any point of intersection $\bar x \in Q \cap M$ satisfying the transversality condition
\bmye \label{transversal manifold}
N_Q(\bar x) \cap -N_{M}(\bar x) ~=~ \{0\},
\emye
where the notation $N_Q(\cdot)$ denotes the (limiting) normal cone to $Q$.  (We refer the reader to \cite{VA} for standard terminology in variational analysis.)  Then, as shown in \cite{alt_cur}, providing the initial point $z^0 \in Q$ is sufficiently close to $\bar x$, the generated sequence of points $z^0,z^1,z^2,\ldots$ in $Q$ converges linearly to a point $\hat z \in Q \cap M$:  there exist positive constants $c < 1$ (a ``rate'' determined by a quantitative version of the transversality condition (\ref{transversal manifold})) and $\rho>0$ such that $|z^k - \hat z| < \rho c^k$ for all iterations $k=0,1,2,\ldots$.  In fact just a weaker notion of {\em intrinsic transversality} \cite{alt_cur} suffices.

This convergence result is striking in its generality, assuming nothing about the closed sets $Q$ and $M$ beyond their transversal intersection.  This underlying assumption is subtle, but in fact a  more intuitive assumption suffices, under a modest geometric assumption about one of the sets.  

Specifically, consider the following idea from \cite{lewis-luke-malick}.

\begin{defn} \label{super-regular}
{\rm
A set $Q \subset {\mathbf X}$ is {\em super-regular} at a point $\bar x$ if, given any angle $\gamma > 0$, for any  points $z \in Q$ and $x \not\in Q$ sufficiently close to $\bar x$, and any point $z' \in P_Q(x)$ with 
$z' \ne z$, the angle between the vectors $z-z'$ and $x-z'$ is at least $\frac{\pi}{2} - \gamma$.
}
\end{defn}

\noindent
In particular, closed convex sets are everywhere super-regular, as are manifolds defined by $C^{1}$-smooth equations with linearly independent gradients.  

Somewhat stronger but more standard than super-regularity is the notion of prox-regularity:  the set $Q$ is {\em prox-regular} at a point $\bx \in Q$ f the projection $P_Q$ is single-valued near $\bx$ (see \cite{loc_diff}).  In that case $P_Q$ is also Lipschitz near $\bx$, and satisfies the following characterization:
\bmye \label{projections}
P_Q(z) = x \quad \Leftrightarrow \quad z-x \in N_Q(x)
\emye
for all points $z,x \in \X$ near $\bx$.  Furthermore, $Q$ is then super-regular at $\bx$, as shown in \cite{lewis-luke-malick}.

Now consider the following simple geometric condition (called ``$0$-separability'' in \cite{noll-rondepierre}), specifically designed with alternating projections in mind.

\begin{defn} \label{separable}
{\rm
A set $Q \subset {\mathbf X}$ intersects a set $M \subset {\mathbf X}$ {\em separably} at a point $\bar x \in M \cap Q$ if there exists an angle $\alpha > 0$ such that, for any point $z \in Q \setminus M$ sufficiently close to 
$\bar x$, and any points $x \in P_M(z) \setminus Q$ and $z' \in P_Q(x)$, the angle between the vectors 
$z-x$ and $z'-x$ is at least $\alpha$.
}
\end{defn}

\noindent
It is easy to check that this property follows from intrinsic transversality, and hence in particular from transversality.  Notice that, unlike those conditions, the separable intersection property is not symmetric in the sets $M$ and $Q$.  

These properties together imply linear convergence of alternating projections, beginning close to the intersection.  To see this, suppose that both Definition \ref{super-regular} and Definition \ref{separable} hold: the set $Q$ intersects the set $M$ separably at the point $\bx$ and is super-regular there.  Now recall the following simple geometric argument, identical to the approach in \cite{lewis-luke-malick}, and extended further in \cite{noll-rondepierre}.  Consider one iteration of the algorithm, starting from a point $z$, and generating the points $x$ and $z'$ in turn.  Given any angle $\gamma > 0$, the triangle near the intersection point 
$\bar x$ defined by the three points $z$, $x$ and $z'$ has angle at least $\alpha$ at $x$ (by separability) and at least $\frac{\pi}{2} - \gamma$ at $z'$ (by super-regularity), and hence at most 
$\frac{\pi}{2} + \gamma - \alpha$ at $z$.  The Sine Rule shows
\[
\frac{|z'-x|}{\sin(\frac{\pi}{2} + \gamma - \alpha)} ~\le~
\frac{|z-x|}{\sin(\frac{\pi}{2} - \gamma)}.
\]
Now fix any constant $\tau$ in the interval $(\cos\alpha,1)$.  By fixing $\gamma$ sufficiently small, we deduce that the next point generated by the algorithm, $x' \in P_M(z')$, satisfies
\bmye \label{reduction}
|z' - x'| ~\le~ |z'-x| ~\le~ \tau |z-x|.
\emye
An easy induction argument now shows linear convergence of the iterate sequence $z_0,z_1,z_2,\ldots$ to a point in the intersection $Q \cap M$.  The rate is $\sqrt{\tau}$, which we can make arbitrarily close to 
$\sqrt{\cos\alpha}$.

\section{Inexact projections}
In general, computing exact projections on nonconvex sets may be an unrealistic subproblem.  In this work we consider to what extent we can relax the projection operations while still retaining the geometry underlying the simple linear convergence argument.  In our approach, we distinguish the roles of the sets $M$ and $Q$.  For simplicity, we assume that we can compute exact nearest points in the set $Q$, which might, for example, be a semidefinite-representable convex set.  On the other hand, we rely only on relaxed projections onto the set $M$, projecting instead onto approximations of $M$ using information from previous iterates.  Earlier work on inexact alternating projection schemes include \cite[Thm 6.1]{lewis-luke-malick} and \cite{kruger-thao}.  In our current development, a direct approach is more convenient.  

We first consider the following generalization of alternating projections, which allows small errors in the computation of the projection $P_M(z)$.

\begin{alg}[Inexact alternating projections]
\label{approximate alternating projections 2}
{\rm
\begin{algorithmic}
\STATE
\STATE  choose $z \in Q$ and $\epsilon \ge 0$; 
\WHILE{$z \not\in M$}
\STATE  choose $x \in {\mathbf X}$ with $d_{P_M(z)}(x) \le \epsilon d_M(z)$;
\STATE  choose $z_+ \in P_{Q}(x)$;
\STATE  $z=z_+$;
\ENDWHILE
\end{algorithmic}
}
\end{alg}

\noindent
When the tolerance $\epsilon$ is zero, we recover the exact alternating projections scheme.  When 
$\epsilon > 0$, however, we no longer even require the point $x$ to lie in the set $M$.  

Suppose that $Q$ intersects $M$ separably at a point $\bar x$, as holds in particular under the usual transversality assumption (\ref{transversal manifold}).  Suppose furthermore that $Q$ is prox-regular (and hence super-regular) at $\bar x$.  Then, providing $\epsilon$ is sufficiently small, starting this scheme at any point $z \in Q$ near $\bar x$ again results in linear convergence to a point in the intersection $Q \cap M$.  To see this, notice that at each iteration there exists a point $x' \in P_M(z)$ satisfying 
$|x-x'| \le \epsilon |z-x'|$.  Denote by $x'_+ \in P_M(z_+)$ a corresponding such point at the next iteration, and let $z' = P_Q(x')$.  As before, after initially fixing any constant $\tau$ in the interval 
$(\cos\alpha,1)$, we deduce $|z' - x'| \le \tau |z-x'|$.  Furthermore, if we fix a local Lipschitz constant $\kappa > 0$ for the projection $P_Q$, then we have 
\[
|z_+ - z'| = |P_Q(x)-P_Q(x')|  \le \kappa|x-x'|  \le \kappa\epsilon|z-x'|.  
\]
An immediate consequence is
\[
|z_+ - x'_+| \le |z_+ - x'| \le |z_+ -z'| + |z'-x'| \le (\kappa\epsilon + \tau)|z-x'|.
\]
Again, providing $\epsilon < \frac{1-\tau}{\kappa}$, we deduce linear convergence of the iterates $z$ to a point in the intersection $Q \cap M$.  As we reduce $\epsilon$ to zero, the convergence rate of the approximate algorithm grows arbitrarily close to that of the exact alternating projection scheme.

With this convergence result in mind, we make the following definition.

\begin{defn} \label{inexact-def}
{\rm
A map $\Phi \colon {\mathbf X} \to {\mathbf X}$ is an {\em inexact projection} onto a closed set $M \subset {\mathbf X}$ around a point $\bar x \in M$ if 
\[
d_{P_M(z)}\big(\Phi(z)\big) = o\big(d_M(z)\big)
\]
as  $z \to \bar x$.
}
\end{defn}

\noindent
We have then proved the following result.

\begin{thm} \label{local-general2}
Suppose a closed set $Q$ intersects a closed set $M$ separably at a point $\bx$, and is prox-regular there.  If $\Phi$ is an inexact projection onto $M$ around $\bx$, then starting from any nearby point $z \in Q$, the inexact alternating projection iteration
\[
z \leftarrow P_Q\big(\Phi(z)\big)
\]
converges linearly to a point in the intersection $M \cap Q$.
\end{thm} 

\section{Inexact projection onto constraint sets}
In this section we show how to compute computationally-cheap inexact projections for constraint sets of the form
\bmye \label{constraint-set}
M ~=~ \{ x \in \X : G(x) \le 0,~P(x) \le 0,~ H(x)=0 \},
\emye
for maps $G \colon \X \to \Rn$, $G \colon \X \to \R^m$ and $H \colon \X \to \Y$, where $\Y$ is a second Euclidean space. We focus on a neighborhood of a point $\bx \in M$ where $P(\bx) < 0$:  nearby, the map 
$P$ models the inactive constraints.  The assumption that the maps $G$, $P$, and $H$ are everywhere defined is just for notational convenience. We make the following classical assumption about the derivative maps 
$\nabla G(\bx) \colon \X \to \Rn$ and $\nabla H(\bx) \colon \X \to \Y$ and their adjoints 
$\nabla G(\bx)^* \colon \Rn \to \X$ and $\nabla H(\bx)^* \colon \Y \to \X$.

\begin{ass}[Linear independence constraint qualification] \label{licq}
{\rm 
\mbox{}
\\
The maps $G$, $P$ and $H$ are $C^{2}$-smooth around a point $\bx \in \X$ satisfying $G(\bx) = 0$, $P(\bx) < 0$, and $H(\bx) = 0$, and they satisfy the condition
\[
\nabla G(\bx)^* w + \nabla H(\bx)^* y = 0 \quad \Rightarrow \quad w=0 ~\mbox{and}~ y=0.
\]
}
\end{ass}  

\noindent
Under these assumptions, the set $M$ is strongly amenable at $\bx$, in the sense of \cite{amen}, and hence in particular prox-regular there \cite[Proposition 2.5]{prox_reg}.  

We can construct an inexact projection for the constraint set $M$ as follows.  For any point $z$ in $\X$ near $\bx$, we consider the projection of $z$ onto the polyhedral approximation to the set $M$ obtained by linearizing its defining constraints at $z$:  in other words, we solve the convex quadratic program
\bmye \label{linearized-constraints}
\left\{
\begin{array}{lrcl}
\mbox{minimize}		& \frac{1}{2}|x-z|^2		&		&	\\
\mbox{subject to}	& G(z) + \nabla G(z)(x-z)	& \le	& 0	\\
					& P(z) + \nabla P(z)(x-z)	& \le	& 0	\\
					& H(z) + \nabla H(z)(x-z)	& =		& 0 \\
					& x							& \in	& \X.
\end{array}
\right.
\emye

\noindent
As a consequence we deduce, under reasonable conditions, local linear convergence of the first algorithm we presented in the introduction.

Our aim, then, is to prove the following result.

\begin{thm}
Consider the constraint set $M$ defined by equation (\ref{constraint-set}).  Suppose the linear independence constraint qualification (Assumption \ref{licq}) holds.  For points $z \in \X$ near the point $\bx$, denote by $\Phi(z)$ the projection onto the linearized constraints, so $\Phi(z)$ is the unique optimal solution of the quadratic program (\ref{linearized-constraints}).  Then, the map $\Phi$ is an inexact projection onto $M$ around $\bx$:
\[
|\Phi(z) - P_M(z)| = O(d_M(z)^2) \quad \mbox{as}~ z \to \bx.
\]
\end{thm}

\pf
Near $\bx$, the constraint involving the continuous map $P$ has no impact on the set $M$, so for the time being we ignore it.   Denote the nonnegative orthant in $\Rn$ by $\Rn_+$, and the set of complementary pairs of vectors by
\[
\Omega ~=~ \{ (w,s) \in \Rn_+ \times \Rn_+ : \ip{w}{s} = 0 \}.
\]
Then any point $x \in \X$ near $\bx$ and vector $p \in \X$ satisfies $x \in M$ and $p \in N_M(x)$ if and only if there exists a pair $(w,s) \in \Omega$ and a vector $y \in \Y$  satisfying $G(x)+s=0$ and $p = \nabla G(x)^* w + \nabla H(x)^* y$.

For any nearby point $z \in \X$, we can now characterize the unique projection $P_M(z)$ using the property (\ref{projections}).  Consider the following set-valued mapping 
\[
\Psi \colon \X \times \Rn \times \Rn \times \Y \tto \Rn \times \Y \times \X  
\]
defined by setting $\Psi(x,w,s,y)$ to be the single point
\[
\left( 
\begin{array}{c}
G(x)+s \\
H(x) \\
x + \nabla G(x)^* w + \nabla H(x)^* y 
\end{array}
\right)^T
\]
if $(w,s) \in \Omega$, and to be empty otherwise.  
Then, for all points $z$ and $x$ near the point $\bx$, the condition
\[
\left(
\begin{array}{c}
0\\
0\\
z
\end{array}
\right)^T 
\in \Psi(x,w,s,y) \quad \mbox{for some}~ w \in \Rn,~ s \in \Rn,~ y \in \Y
\]
is necessary and sufficient for $x = P_M(z)$. 

We claim that the map $\Psi$ is metrically regular at $(\bx,0,0,0)$ for $(0,0,\bx)$.  To see this, we linearize around the point $\bx$ to obtain the following new set-valued mapping 
$\hat \Psi \colon \X \times \Rn \times \Rn \times \Y \tto \Rn \times \Y \times \X$.  
Define the set
$\hat \Psi(x,w,s,y)$ to be the single point
\[
\left(
\begin{array}{c}
\nabla G(\bx)(x-\bx) + s\\
\nabla H(\bx)(x-\bx)\\
x + \nabla G(\bx)^* w + \nabla H(\bx)^* y
\end{array}
\right)^T
\]
if $(w,s) \in \Omega$, and to be empty otherwise.  This map differs from the map $\Psi$ by the map that takes $(x,w,s,y) \in \X \times \Rn \times \Rn \times \Y$ to the point
\[
\left(
\begin{array}{c}
G(x) - \nabla G(\bx)(x-\bx)\\
H(x) - \nabla H(\bx)(x-\bx)\\ 
\big(\nabla G(x) - \nabla G(\bx)\big)^* w + \big(\nabla H(x) - \nabla H(\bx)\big)^* y
\end{array}
\right)^T.
\]
This latter map is continuously differentiable around the point $(\bx,0,0,0)$, with zero derivative.  Hence, using the coderivative criterion \cite{VA}, metric regularity of the maps $\Psi$ and $\hat \Psi$ at $(\bx,0,0,0)$ for for $(0,0,\bx)$ are equivalent.

For all $(u,v,z) \in \Rn \times \Y \times \X$ and $x \in \X$, the condition
\[
\left(
\begin{array}{c}
u\\
v\\
z
\end{array}
\right)^T 
\in \hat\Psi(x,w,s,y) \quad \mbox{for some}~ w \in \Rn,~ s \in \Rn,~ y \in \Y
\]
is just the necessary and sufficient first-order condition for $x$ to be an optimal solution for the convex quadratic program
\[
\left\{
\begin{array}{lrcl}
\mbox{minimize}		& \frac{1}{2}|x-z|^2	&		&	\\
\mbox{subject to}	& \nabla G(\bx)x		& \le	& \nabla G(\bx)\bx + u	\\
					& \nabla H(\bx)x		& =		& \nabla H(\bx)\bx + v \\
					& x						& \in	& \X.
\end{array}
\right.
\]
This problem is feasible, by the linear independence constraint qualification, and hence has a unique optimal solution $\xi(u,v,z)$, which satisfies the first-order condition.  Hence the set 
$\hat\Psi^{-1}(u,v,z)$ is nonempty, and every element $(x,w,s,y) \in \hat\Psi^{-1}(u,v,z)$ must satisfy
\begin{eqnarray*}
x &=& \xi(u,v,z)\\
s &=& u - \nabla G(\bx)(x-\bx)\\ 
\nabla G(\bx)^* w + \nabla H(\bx)^* y &=& z - x.
\end{eqnarray*}
Thus the point $z-\xi(u,v,z)$ must lie in the range 
${\mathbf S} = \nabla G(\bx)^* \Rn + \nabla H(\bx)^* \Y$.  By the linear independence constraint qualification, the linear map 
\[
(w,y) \mapsto \nabla G(\bx)^* w + \nabla H(\bx)^* y
\]
has a linear inverse $K \colon {\mathbf S} \to \Rn \times \Y$.  We deduce that $\hat\Psi^{-1}(u,v,z)$ is a singleton $\{(x,w,s,y)\}$, where
\begin{eqnarray*}
x &=& \xi(u,v,z)\\
s &=& u - \nabla G(\bx)(\xi(u,v,z)-\bx)\\
(w,y) &=& K(z-\xi(u,v,z)).
\end{eqnarray*}

For any fixed point $z \in X$, the dependence of the optimal solution $\xi(u,v,z)$ on the parameters 
$(u,v) \in \Rn \times \Y$ is globally Lipschitz, with constant independent of $z$ (see \cite{yen95}, and also the discussion in Facchinei-Pang, Section 4.7), and on the other hand, for fixed $(u,v)$, as a function of the point $z$, it is globally $1$-Lipschitz (being a projection).  Hence the function $\xi$ is globally Lipschitz.  Consequently the map $\hat\Psi^{-1}$ is single-valued and Lipschitz, so the map $\hat\Psi$ is metrically regular, and hence so is the map $\Psi$ at $(\bx,0,0,0)$ for for $(0,0,\bx)$, as we claimed.

Equipped with this metric regularity property, we can now prove the desired inexact projection property for the constraint set $M$ defined by equation (\ref{constraint-set}).  For any point $z$ in $\X$ near 
$\bx$, we solve the convex quadratic program (\ref{linearized-constraints}).
To see that this calculation results in an inexact projection, in the sense of Definition~\ref{inexact-def}, we can ignore the inactive constraints modeled by the map $P$, and consider the unique optimal solution $x_z$ of the convex quadratic program
\[
\left\{
\begin{array}{lrcl}
\mbox{minimize}		& \frac{1}{2}|x-z|^2		&		&	\\
\mbox{subject to}	& G(z) + \nabla G(z)(x-z)	& \le	& 0	\\
					& H(z) + \nabla H(z)(x-z)	& =		& 0 \\
					& x							& \in	& \X.
\end{array}
\right.
\]
Define a $C^{2}$-smooth map $A \colon \X \to \Rn \times \Y$ by $A(z) = \big(G(z),H(z)\big)$, for points $z \in \X$.  By the linear independence constraint qualification, the linear map 
$\nabla A(z) \colon \X \to \Rn \times \Y$ is one-to-one for all $z$ near $\bx$, and hence the composition
\[
\nabla A(z)\nabla A(z)^* \colon \Rn \times \Y \to \Rn \times \Y
\]
is positive definite.  The point
\[
x = z - \nabla A(z)^* \big(\nabla A(z)\nabla A(z)^*\big)^{-1}A(z)
\] 
is feasible for the quadratic program.  Hence we deduce
\[
|x_z - z| ~\le~ \Big| \nabla A(z)^* \big(\nabla A(z)\nabla A(z)^*\big)^{-1}A(z) \Big| ~\rightarrow~ 0
\]
as $z \to \bx$, so $x_z \to \bx$.  Furthermore, for all $z$ near $\bar x$, the optimal solution $x_z$ is also feasible for the the linearized problem (\ref{linearized-constraints}), so it is also the unique optimal solution of that problem.

Using the first-order optimality conditions, there exist a unique slack vector $s_z \in \Rn$ and multiplier vectors 
$w_z \in \Rn$ and $y_z \in \Y$ such that $(w_z,s_z) \in \Omega$ and 
\begin{eqnarray*}
G(z) + \nabla G(z)(x_z-z) + s_z				& = & 0 \\
H(z) + \nabla H(z)(x_z-z)					& = & 0 \\
x_z + \nabla G(z)^*w_z + \nabla H(z)^*y_z	& = & z.
\end{eqnarray*}
In fact $w_z$ and $y_z$ are also unique, by the linear independence constraint qualification, which furthermore implies
\bmye \label{order-one}
|w_z| + |y_z| = O(|x_z-z|)
\emye
as $z \to \bx$.  The $C^{2}$-smoothness of $G$ and $H$ implies
\begin{eqnarray*}
G(x_z) +  s_z	& = & O(|x_z-z|^2) \\
H(x_z) 			& = & O(|x_z-z|^2).
\end{eqnarray*}
Furthermore we have
\begin{eqnarray*}
\lefteqn{x_z + \nabla G(x_z)^*w_z + \nabla H(x_z)^*y_z - z } \hspace{2cm} & & \\
& = & (\nabla G(x_z) - \nabla G(z))^*w_z + (\nabla H(x_z) - \nabla H(z))^*y_z \\
& = & O(|x_z-z|) |w_z| + O(|x_z-z|) |y_z| \\
& = & O(|x_z-z|^2)
\end{eqnarray*}
by equation (\ref{order-one}).  By the metric regularity of the map $\Psi$, we deduce 
\[
|x_z - P_M(z)| = O(|x_z-z|^2) \quad \mbox{as}~ z \to \bx.
\]
Notice the inequalities
\begin{eqnarray*}
\lefteqn{|x_z-z| \le |x_z - P_M(z)| + |P_M(z) - z| } \hspace{2cm} & & \\
& = & O(|x_z-z|^2) + d_M(z) \le \frac{1}{2}|x_z-z| + d_M(z)
\end{eqnarray*}
for all $z$ near $\bx$, which imply $|x_z-z| \le 2d_M(z)$, and hence
\[
|x_z - P_M(z)| = O(d_M(z)^2) \quad \mbox{as}~ z \to \bx.
\]
This completes the proof that the map $z \mapsto x_z$ is an inexact projection. 
\finpf 
%
%

When the set $M$ is defined only by equalities, the assumptions imply that it is a manifold around the point $\bx$.  As we note in the next section, manifolds also have dual representations, in terms of local coordinates.  Given such a representation, as we see next, another approach to approximate alternating projections is available, which in particular allows us to dispense with the assumption that the set $Q$ is prox-regular.

\section{Approximate projections:  a second approach}
We return to the basic alternating projections convergence scheme that we outlined in Section~\ref{alternating section}.  In contrast with the previous section, we consider a different approach to errors in the computation of the projection $P_{M}(y)$.  Shifting notation slightly, we again consider alternating projections for two closed sets $M$ and $Q$, this time in a Euclidean space $\Y$.  We suppose, at a common point $\bar y$, that $M$ intersects $Q$ separably (as opposed to vice versa, as before) and that $M$ is super-regular there (as opposed to prox-regularity of 
$Q$, as before).

Consider the following algorithm.

\begin{alg}[Approximate alternating projections]
\label{approximate alternating projections}
{\rm
\begin{algorithmic}
\STATE
\STATE  choose $z \in M$ and $\epsilon \ge 0$; 
\WHILE{$z \not\in Q$}
\STATE  choose $y \in P_Q(z)$;
\STATE	choose $z' \in M$ with $d_{P_M(y)}(z') \le \epsilon |y-z|$;
\STATE  $z=z'$;
\ENDWHILE
\end{algorithmic}
}
\end{alg}

\noindent
Notice that, unlike in the previous inexact alternating projection scheme (Algorithm~\ref{approximate alternating projections 2}), we now require the approximate projection $z'$ to lie in the set $M$.

Our assumptions again imply linear convergence of this algorithm, beginning close to the intersection point 
$\bar y$, providing that the tolerance $\epsilon$ is sufficiently small.  To see this, as before consider one iteration of the algorithm, starting from a point $z$, and generating the points $y$ and $z'$ in turn, choose $z'' \in P_{P_M(y)}(z')$
and let $y' \in P_Q(z')$ denote the next point generated.  For any constant $\tau$ in the interval 
$(\cos\alpha,1)$, the analog of inequalities (\ref{reduction}) still hold,
\[
|z'' - y| \le \tau|z-y|,
\]
and so we have
\[
|z'-y'| \le |z' - y| \le |z' - z''| + |z'' - y| \le (\epsilon + \tau)|z-y|.
\]
For any nonnegative $\epsilon < 1-\tau$, we deduce linear convergence of the iterate sequence of iterates $z$ to a point in the intersection $Q \cap M$, now with rate $\sqrt{\tau + \epsilon}$, which again we can make arbitrarily close to $\sqrt{\cos\alpha}$.

With this convergence result in mind, we make the following definition, capturing the idea of approximating the projection of a point $y \in \Y$ onto a set $M \subset Y$ by an inexact estimate 
$\Phi(z,y)$ that depends on a base point $z \in M$.

\begin{defn}
{\rm
A {\em faithful approximation} of the projection onto a set $M \subset \Y$ around a point $\bar y \in M$ is a map $\Phi \colon V \times \Y \to M$, where $V$ is a neighborhood of $\bar y$ in $M$, with the following property:  given any constant $\epsilon > 0$ and angle $\alpha > 0$, if
points $z \in M$ and $y \not\in M$ are sufficiently close to $\bar y$, then any point $\hat z \in P_M(y)$ such that the angle between the vectors $z-y$ and $\hat z-y$ is at least $\alpha$ must satisfy 
$|\hat z-\Phi(z,y)| \le \epsilon |y-z|$.
}
\end{defn}

\noindent
We have then proved the following result.

\begin{thm} \label{local-general}
Suppose a closed set $M$ intersects a closed set $Q$ separably at a point $\bar y$, and is super-regular there.  If $\Phi$ is a faithful approximation of the projection onto $M$ around $\bar y$, then starting from any nearby point $z \in M$, the 
approximate alternating projection iteration
\[
z \leftarrow \Phi(z,y) \quad \mbox{for any}~ y \in P_Q(z)
\]
converges linearly to a point in the intersection $M \cap Q$.
\end{thm}

\section{Smooth inclusions}
We apply the approximate alternating projections framework of the previous section to the following problem.  Given two Euclidean spaces $\X$ and $\Y$, a continuously differentiable map 
$F \colon \X \to \Y$, and a closed set 
$Q \subset \Y$, we consider the problem of finding a point $x \in \X$ satisfying
\[
F(x) \in Q.
\]
We can then consider the following simple splitting algorithm, which alternates between projections onto the set $Q$ and linear least squares problems.

\begin{alg}[Linearized alternating projections]
\label{inexact projections}
{\rm
\begin{algorithmic}
\STATE
\STATE  choose $x \in \X$; 
\WHILE{$F(x) \not\in Q$}
\STATE	choose $y \in P_Q\big(F(x)\big)$;
\STATE  choose $s \in \X$ minimizing $|F(x) + \nabla F(x)s - y|$;
\STATE	$x=x+s$;
\ENDWHILE
\end{algorithmic}
}
\end{alg}

\noindent
In the case $Q = \{0\}$, this is simply the Gauss-Newton method for the nonlinear least squares problem of minimizing $|F(\cdot)|^2$.  On the other hand, if we replace the linear approximation 
$F(x) + \nabla F(x)s$ in the algorithm by the exact value \mbox{$F(x+s)$}, then, under reasonable conditions, the method becomes the exact alternating projection algorithm of Section \ref{alternating section} for the intersection of $Q$ with a certain smooth manifold, namely the image under the map $F$ of an open set containing the iterates $x$.  We can then deduce linear convergence close to well-behaved solutions.  Here, we parallel that development to prove local linear convergence of Algorithm~\ref{inexact projections} by interpreting it as a version of the approximate projection scheme in the previous section:  Algorithm \ref{approximate alternating projections}.

More formally, we summarize the assumptions we need below.
\begin{ass} \label{assumptions}
{\rm
The map $F$ is continuously differentiable around the point $\bx$, with one-to-one derivative 
$\nabla F(\bx)$, and satisfies $F(\bar x) \in Q$.  Furthermore, there exists an angle 
$\alpha > 0$ such that whenever two vectors $x,x' \in \X$ near $\bar x$ satisfy
\[
\nabla F(x')^*(y - F(x')) = 0 \quad \mbox{and} \quad F(x) \ne y \ne F(x') \quad 
\mbox{for some} \quad y \in P_Q\big(F(x)\big),
\]
then the angle between the vectors $F(x)-y$ and $F(x')-y$ is at least $\alpha$.
}
\end{ass}

Under this assumption, we study Algorithm \ref{inexact projections} in a neighborhood of a particular solution $\bar x$ of our problem.  The one-to-one assumption on the derivative in particular allows us to work with the following basic terminology for manifolds.

\begin{defn} \label{manifold}
{\rm 
A set $M \subset \Y$ is a {\em manifold} around a point $\bar y \in M$ when there exists a corresponding {\em local coordinate system\/}: a continuously differentiable map $F$ from another Euclidean space $\X$ into $\Y$ whose restriction to some open set $U \subset \X$ is a diffeomorphism onto an open neighborhood $V$ of $\bar y$ in $\cM$.
}
\end{defn}

\noindent
Assumption \ref{assumptions} hence guarantees that, for any sufficiently small neighborhood 
$U$ of $\bx$, the image $M = F(U)$ is a manifold around the point $\bar y = F(\bx)$, and the map $F$ is a local coordinate system.  Furthermore, the normal space to the manifold is given by
\[
N_M\big(F(x)\big) = \mbox{Null}\big( \nabla F(x)^* \big),
\]
for all $x$ near $\bx$.
Thus the angle condition is exactly the property that $M$ intersects the set $Q$ separably at the point 
$\bar y$.  It holds in particular under the usual transversality condition (\ref{transversal manifold}), which here becomes
\[
N_Q\big(F(\bx)\big) \cap N\big(\nabla F(\bx)^*\big) = \{0\}.
\]

Our aim is to prove the following result.

\begin{thm}[Local linear convergence] \label{local}
If Assumption \ref{assumptions} holds, then, starting from any point $x \in \X$ sufficiently close to the solution $\bar x$, the linearized alternating projections method, Algorithm \ref{inexact projections}, 
converges linearly to a point in the solution set $F^{-1}(Q)$.
\end{thm} 

\section{Faithful approximate projections on manifolds}
As in Definition \ref{manifold}, consider a local coordinate system $F$ for a manifold $M$ around a point 
$\bar y \in M$.  We can calculate approximate projections onto the manifold $M$ via a map $\Phi \colon V \times \Y \to \Y$ defined as follows.  For coordinate vectors $x \in U$ and points $y \in \Y$, define
\bmye \label{approx}
\Phi(F(x),y) = F(x+s) \quad \mbox{where $s \in \Y$ minimizes} \quad |F(x) + \nabla F(x) s - y|.
\emye
Notice that the vector $s$ is unique, since the derivative map $\nabla F(x)$ is one-to-one.  Calculating $s$ is a simple least squares problem.  If $y$ is close to the point $F(x)$, then the vector $s$ is small, and hence the point $F(x+s)$ lies on the manifold $M$, and is close to the point 
$F(x) + \nabla F(x) s$, which is the projection of $y$ onto the tangent space $T_{M}(x)$. In this sense, the map $\Phi(F(x),\cdot)$ approximates the projection mapping $P_{M}$ near the point $F(x)$.  We can make this statement more precise as follows.

\begin{thm}
The map $\Phi$ defined by equation (\ref{approx}) is a faithful approximation of the projection for the manifold $M$ around the point $\bar y$.
\end{thm}

\pf
Suppose the result fails, so there exists an angle $\alpha > 0$ and a constant $\epsilon > 0$, and sequences
$z_k \in M$, $y_k \in \Y$, and $z_k' \in P_M(y_k)$, all converging to $\bar y$, and satisfying
\[
\epsilon_k = |z_k'-z_k| \ge |y_k-z_k|\sin\alpha \quad \mbox{and} \quad
|z_k'-\Phi(z_k,y_k)| > \epsilon |y_k-z_k|.
\]
Note the inequality
\bmye \label{triangle}
\epsilon_k \le |z_k' - y_k| + |y_k - z_k| \le 2|y_k - z_k|.
\emye

Denote the coordinates of $z_k,z'_k$ by $x_k,x'_k \in U$ respectively (clearly distinct), so 
$F(x_k) = z_k$ and $F(x'_k) = z'_k$.   Since $F$ is a diffeomorphism between $U$ and $V$, the sequences 
$(x_k)$ and $(x'_k)$ both converge to the common limit $\bx \in U$, the coordinates of $\bar y$.  Define $\delta_k = |x'_k - x_k|$:  the diffeomorphism property ensures $\epsilon_k = O(\delta_k)$ and $\delta_k = O(\epsilon_k)$.  After taking a subsequence, we can suppose the unit vectors $\frac{1}{\delta_k} (x'_k - x_k)$ converge to a unit vector $u$. 

Now define vectors $w_k = y_k-z'_k \in  N_{M}(z'_k)$.  Notice $\nabla F(x'_k)^* w_k = 0$, and
\[
|w_k| \le |y_k-z_k| + |z_k - z'_k| = O(\epsilon_k) + \epsilon_k = O(\epsilon_k) = O(\delta_k).
\]
Taking a further subsequence, we can then suppose that the vectors $\frac{1}{\delta_k} w_k$ converge to some vector in the null space of $\nabla F(\bx)^*$.

By definition, $\Phi(z_k,y_k) = F(x_k + s_k)$, where the vector $s_k \in \X$ minimizes the function
\[
s \mapsto |F(x_k) + \nabla F(x_k) s - y_k|.
\]
We deduce
\[
\nabla F(x_k)^*(F(x_k) + \nabla F(x_k) s_k - F(x'_k) - w_k) = 0,
\]
so
\[
s_k = \big(\nabla F(x_k)^* \nabla F(x_k)\big)^{-1}\nabla F(x_k)^* \big(F(x'_k) - F(x_k) + w_k\big).
\]
We deduce $\frac{1}{\delta_k} s_k \to u$.  Hence
\[
(x_k + s_k) - x'_k = o(\delta_k)
\]
so inequality (\ref{triangle}) implies
\[
\Phi(z_k,y_k) - z_k' = F(x_k + s_k) - F(x'_k) = o(\delta_k) = o(\epsilon_k) = o(|y_k - z_k|)
\]
which is a contradiction.
\finpf

\noindent
Theorem \ref{local} (Local linear convergence) now follows directly from Theorem \ref{local-general}.

%
%

\def\cprime{$'$} \def\cprime{$'$}

\end{document}